\documentclass[a4paper,12pt]{article}
\usepackage[english]{babel}
\usepackage[T2A]{fontenc}
\usepackage[cp1251]{inputenc}
\usepackage[english]{babel}
\usepackage{amsthm}
\usepackage[tbtags]{amsmath}
\usepackage{amsfonts,amssymb}
\begin{document}

\begin{center}
{\large\bf On Rings over Which All Modules Are Direct Sums of Distributive Modules}
\end{center}
\begin{center}
A.A. Tuganbaev\footnote{National Research University MPEI, Moscow, Lomonosov Moscow State University; email: tuganbaev@gmail.com .}
\end{center}

\textbf{Abstract.} A module is said to be \textsf{distributive} if the lattice of its submodules is distributive. A direct sum of distributive modules is called a \textsf{semidistributive} module. In this paper we consider rings $A$ such that all right $A$-modules are semidistributive.


The work is supported by Russian Scientific Foundation, project 22-11-00052.

\textbf{MSC2020 database 16D25, 16R99}

\section{Introduction}\label{sec1}

We only consider associative unital non-zero rings and unitary modules. The words of type <<a right Artinian ring $A$>> (correspondingly, <<an Artinian ring $A$>>) mean that the module $A_A$ is Artinian (correspondingly, the both modules $A_A$ and $_AA$ are Artinian).

In the book <<Dniester Notebook>> \cite{Dni93} consisting of some problems of ring theory, L.A.~Skornyakov posed the following two questions (Problem 1.116): Over which rings all right modules are semidistributive? Are there non-Artinian rings with this property? 

\textbf{Remark 1.1.} In \cite{Tug80}, it is proved that a ring, over which all right modules are semidistributive, is a right Artinian, right K\"othe ring. (A ring $A$ is called a \textsf{right K\"othe ring} if every right $A$-module is a direct sum of cyclic modules.) It is well-known that right K\"othe rings are Artinian. This gives negative answer to the second question of L.A.~Skornyakov; also see \cite[Theorem 11.6]{Tug98}. Therefore, all rings, over which all right modules are semidistributive, are right K\"othe, Artinian rings. 
The problem of describing the right K\"ete rings $A$, restricting ourselves only to the internal properties of the ring $A$, remains unsolved for arbitrary rings; this problem is called the \textsf{K\"othe's problem}.

\textbf{Remark 1.2.} See \cite{Tug80} and \cite[Section 11.1]{Tug98} on the second question of L.A.~Skornyakov. In \cite{Ful78}, it was studied a partial case of rings, over which all right modules semidistributive.

For a module $M$, we denote by $J(M)$ the Jacobson radical of the module $M$. Let $A$ be a semiprimary ring\footnote{A ring $A$ is said to be \textsf{semiprimary} if the radical $J(A)$ is nilpotent and the factor ring $A/J(A)$ is a semisimple Artinian ring. Every right or left Artinian ring is semiprimary.} and let $B$ be the \textsf{basic} ring of the semiprimary ring $A$, i.e. $B=eAe$, where $e$ is a \textsf{basic} idempotent of the ring $A$; this means that $e=e_1+\ldots+e_n$, where $\{e_1,\ldots,e_n\}$ is a set of local orthogonal idempotents of $A$ such that $\{e_1A,\ldots,e_nA\}$ (correspondingly, $\{Ae_1,\ldots,Ae_n\}$) is the set of all pair-wise non-isomorphic indecomposable direct summands of the module $A_A$ (correspondingly, $_AA$). If $A=B$, then the semiprimary ring $A$ is said to be \textsf{self-basic}. It is well known that the basic ring $B$ of the semiprimary ring $A$ is self-basic and the category $\text{Mod }A$ of all right $A$-modules is equivalent to the category of all right $B$-modules, i.e. the rings $A$ and $B$ are \textsf{Morita equivalent}. In addition, it is known that the property to be an Artinian self-basic ring is preserved under Morita equivalences.

\section{Main Results}\label{sec2}

\textbf{Remark 2.1.} It follows from Remark 1.1 and the above that any ring $A$, over which all right modules are semidistributive, is an Artinian ring with self-basic basic ring $B$ and the rings $A$ and $B$ are Morita equivalent. It is clear that the semidistributivity property of all right modules is preserved under Morita equivalences. Therefore, all right $B$-modules are semidistributive. Therefore, when studying rings over which all right modules are semidistributive, we can restrict ourself by Artinian self-basic K\"othe rings $B$ such that every ring, which is Morita equivalent to the ring $B$, is an Artinian self-basic ring.

A module is said to be \textsf{completely cyclic} if all its submodules are cyclic. 

\textbf{Theorem 2.2 \cite[Theorem 11.6]{Tug98}.} 
For a ring $A$, the following conditions are equivalent.

\textbf{1)} All right $A$-modules semidistributive.

\textbf{2)} $A$ is an Artinian ring and every right $A$-module is a direct sum of completely cyclic distributive modules with composition series.

\textbf{3)} $A$ is an Artinian ring with basic idempotent $e$ and for any right $A$-module $M$, the right $eAe$-module $Me_{eAe}$ is a direct sum of completely cyclic modules with composition series.

A ring $A$ is said to be a ring of \textsf{of finite representation type} if $A$ is an Artinian ring which has up to isomorphism only a finite number of indecomposable right modules and a finite number of indecomposable left modules. 
A ring $A$ is called a \textsf{right Kawada ring} if every ring, which is Morita equivalent to the ring $A$, is a right K\"othe ring. We note that all right K\"othe rings (in particular, all right Kawada rings) are rings of finite representation type.

\textbf{Theorem 2.3 \cite{Tug23}.} For a ring $A$, the following conditions are equivalent.

\textbf{1)} All right $A$-modules semidistributive.

\textbf{2)} $A$ is a right Kawada ring with basic ring $B$ and every right $A$-module and every right $B$-module are direct sums of completely cyclic distributive modules.

\textbf{3)} $A$ is a right Kawada ring with basic ring $B$ and every indecomposable right $B$-module is a completely cyclic module.

\textbf{Remark 2.4.} In \cite{Kaw61}, Kawada solved the K\"othe's problem for finite-dimensional algebras over a field. The Kawada's theorem completely describes self-basic finite-dimensional algebras $A$ over a field such that every indecomposable $A$-module has the square-free socle and the square-free top; in the same work, all indecomposable $A$-modules are described. The K\"othe's problem remains unsolved in general case. 
In \cite{Kaw61}, there are 19 rather complicated conditions for local idempotents of the basic ring $B$ of the algebra $A$; all these conditions hold if and only if $A$ is a right K\"othe ring. In \cite{Rin81}, this Kawada's result is analyzed and commented. To the mentioned 19 conditions, we can add the following condition 20: for any local idempotent $e$ of the basic algebra $B$, the module $eB_B$ is completely cyclic. Therefore, we obtain a formal description of finite-dimensional algebras over a field, over which all right modules are semidistributive. Of coarse, such a description is not very useful.

A ring is said to be \textsf{normal} or \textsf{Abelian} if all its idempotents are central.

\textbf{Corollary 2.5.} If the factor ring $A/J(A)$ of the ring $A$ is normal, then all right $A$-modules are semidistributive if and only if $A$ is an Artinian ring and every right $A$-module is a direct sum of completely cyclic modules with composition series.

A ring, in which all right ideals and all left ideals are principal, is called a \textsf{principal ideal ring}.

\textbf{Theorem 2.6 \cite{Koe35}.} An Artinian principal ideal ring is a K\"othe ring.

\textbf{Corollary 2.7 \cite{CohK51}, \cite{Koe35}.} A commutative ring $A$ is a K\"othe ring if and only if $A$ is an Artinian principal ideal ring.

A module is said to be \textsf{uniserial} if all its submodules are linearly ordered with respect to inclusion. A direct sum of uniserial modules is called a \textsf{serial} module.

With the use of Corollary 2.5 and Theorem 2.6, it is easy to verify Theorem 2.8.

\textbf{Theorem 2.8.} For a normal ring $A$, the following conditions are equivalent.

\textbf{1)} All right $A$-modules are semidistributive.

\textbf{2)} All left $A$-modules are semidistributive.

\textbf{3)} $A$ is an Artinian principal ideal ring.

\textbf{4)} The ring $A$ is isomorphic to the finite direct product of Artinian uniserial rings.

For any module $M$, the \textsf{top} $\text{top }M$ is the factor module $M/J(M)$. A module, which does not contain direct sums of two non-zero isomorphic submodules, is called a \textsf{square-free} module.

\textbf{Theorem 2.9 \cite{BehGMS14}.} A normal ring $A$ is a K\"othe ring if and only if $A$ is an Artinian principal ideal ring.

According to \cite{AsgBK23}, a ring $A$ is called a 
\textsf{strongly right} (correspondingly, \textsf{strongly left}) \textsf{K\"othe ring} if every non-zero right (correspondingly, left) $A$-module is a direct sum of modules with non-zero cyclic square-free top. Right and left strongly K\"othe rings are called \textsf{strongly K\"othe rings}.

According to \cite{AsgBK23}, a ring $A$ is called a 
\textsf{right very strongly} (correspondingly, \textsf{left very strongly}) \textsf{K\"othe ring} if every non-zero right (correspondingly, left) $A$-module is a direct sum of modules with simple top. Right and left very strongly K\"othe rings are called \textsf{very strongly K\"othe rings}.

The following proper inclusions are known (e.g., see \cite{AsgBK23}).
$$
\text{Very strongly right K\"othe rings}\subsetneq
$$
$$
\subsetneq \text{strongly right K\"othe rings}
\subsetneq\text{right K\"othe rings}.
$$

\textbf{Theorem 2.10 \cite{AsgBK23}.} For a ring $A$, the following conditions are equivalent.

\textbf{1)} $A$ is a right K\"othe ring.

\textbf{2)} Every non-zero right $A$-module is a direct sum of modules with non-zero cyclic top.

\textbf{3)} The ring $A$ is right Artinian and every right $A$-module is a direct sum of modules with cyclic top.

\textbf{4)} $A$ is a ring of finite representation type and every (finitely generated) indecomposable right $A$-module has the cyclic top.

\textbf{5)} $A$ is a ring of finite representation type and the top  of every indecomposable right $A$-module $U$ is isomorphically embedded in $A/J$.

\section{Addendum}\label{sec3}

In \cite{AsgBK23b}, the authors defined co-K\"othe rings which are close to K\"othe rings: a ring $A$ is called a \textsf{right} (correspondingly, \textsf{right strongly, right very strongly}) \textsf{co-K\"othe ring} if every non-zero right $A$-module is a direct sum of modules with non-zero cyclic socle (correspondingly, with non-zero square-free socle, with simple socle). The left-side analogues of these notions are defined similarly.

\textbf{Remark 3.1.} By \cite{AsgBK23b}, there is a very strongly right co-K\"othe ring, over which there exists a non-semidistributive right module.

A module is said to be \textsf{uniform} if the intersection of any two its non-zero submodules is not equal to zero.

\textbf{Example 3.2 (see also \cite{Nak40}, \cite{Nak41}, \cite[Example 1.22]{Tug98}, \cite{AsgBK23b}).}\\

Let $A$ be the $5$-dimensional algebra over the field $\mathbb{Z}_2=\mathbb{Z}/2\mathbb{Z}$ consisting of all $3\times 3$-matrices of the form
$$
\begin{bmatrix}
f_{11}&f_{12}&f_{13}\\
0&f_{22}&0\\
0&0&f_{33}
\end{bmatrix},
$$
where $f_{ij}\in \mathbb{Z}_2$. Let $e_{ij}$ be the matrix whose $ij$th entry is equal to $1$ and all other entries
are equal to $0$. Then $\{e_{11}, e_{12}, e_{13}, e_{22}, e_{33}\}$ is a $\mathbb{Z}_2$-basis of the $\mathbb{Z}_2$-algebra $A$. Moreover, the following assertions are true.

\textbf{1.} $1=e_{11} +e_{22} +e_{33}$, where $e_{11}$, $e_{22}$, $e_{33}$ are primitive orthogonal idempotents, $e_{12}\mathbb{Z}_2=e_{12}A$, $e_{13}\mathbb{Z}_2=e_{13}A$, $J(A)=e_{12}\mathbb{Z}_2 + e_{13}\mathbb{Z}_2$, $J^2=0$, and the ring $A/J$ is isomorphic to a direct product of three copies of the field $\mathbb{Z}_2$.

\textbf{2.} $A_A=e_{11}A\oplus e_{22}A\oplus e_{33}A$, where $e_{22}A=e_{22}\mathbb{Z}_2$ and $e_{33}A=e_{33}\mathbb{Z}_2$ are simple projective right $A$-modules which are isomorphic to the modules $e_{12}A $ and $e_{13}A$, respectively.

\textbf{3.} $e_{11}A=e_{11}\mathbb{Z}_2+e_{12}\mathbb{Z}_2+e_{13}\mathbb{Z}_2$ is an indecomposable distributive (and hence square-free) Noetherian Artinian completely cyclic $A$-module, but it is not uniform and every proper non-zero submodule of $e_{11}A$ coincides either with the projective module $e_{12}A\oplus e_{13}A$ or with one of the simple projective non-isomorphic modules $e_{12}A$ and $e_{13}A$.

It can be checked that $A$ is a hereditary Artinian basic ring. So by \cite{Faz20}, $A$ is a K\"othe ring. It can be shown that $A$ has 32 elements that
$|U(A)|=4$ and non-unit elements of $A$ form 3 isomorphism classes of cyclic indecomposable modules, which are as follows:
$$
Q_1=e_{11}A,\quad Q_2=e_{11}A/e_{13}A,\quad Q_3=e_{11}A/e_{12}A,
$$
$$
Q_4=e_{11}A/J(A),\quad Q_5=e_{22}A,\quad Q6=e_{33}A.
$$

Also, every indecomposable cyclic right $A$-module has a square-free socle, since $\text{Soc }(Q_1)=J(A)$ is square-free, $\text{Soc }Q_2 =J(A)e_{13}A$, $\text{Soc }Q3 =J(A)/e_{12}A$, and $Q_4$, $Q_5$ and $Q_6$ are simple.
Then every right $A$-module is a direct sum of square-free modules ($A$ is a (strongly) right co-K\"othe ring) while the right $A$-module $e_{11}A$ is not a direct sum of uniform modules ($A$ is not a very strongly right co-K\"othe ring), since $A$ is not right serial.

We give several results from \cite{AsgBK23b}. 
We recall that a semi-primary ring $A$ is called a \textsf{right QF-2 ring} (resp., \textsf{right co-QF-2 ring}) if every indecomposable projective right $A$-module has a simple essential socle (resp., a simple top). 

A semi-perfect ring $A$ is a generalized left co-QF-2 ring if every indecomposable projective left $A$-module $P$ has a square-free top.
A semi-perfect ring $A$ is a generalized co-QF-2 ring if $A$ is a generalized left and a generalized right co-QF-2 ring.

\textbf{Theorem 3.3 \cite{AsgBK23b}.}\\
The following conditions are equivalent for a ring $A$:
 
\textbf{1)} $A$ is a right co-K\"othe ring.

\textbf{2)} Every non-zero right $A$-module is a direct sum of modules with non-zero top and cyclic essential socle.

\textbf{3)} $A$ is of finite representation type and every (finitely generated) indecomposable right $A$-module has a cyclic (essential) socle.

\textbf{4)} $A$ is of finite representation type and the socle of every indecomposable right $A$-module $U$ is isomorphically embedded in $A/J$.

Let $A$ be a ring of finite representation type and let $U = U_1\oplus\ldots\oplus U_n$ and $\{U_1,\ldots ,U_n\}$ be a complete set of representatives of the isomorphic classes of finitely generated indecomposable right $A$-modules. The \textsf{right Auslander ring} of $A$ is $T=\text{End }U_A$.

\textbf{Theorem 3.4 \cite{AsgBK23b}.}\\
The following conditions are equivalent for a ring $A$ with $J=J(A)$.

\textbf{1)} $A$ is a strongly right co-K\"othe ring.

\textbf{2)} Every right $A$-module is a direct sum of square-free modules.

\textbf{3)} Every non-zero right $A$-module is a direct sum of modules with non-zero top and square-free (cyclic) socle.

\textbf{4)} $A$ is of finite representation type and every (finitely generated) indecomposable right $A$-module has a square-free (cyclic) socle.

\textbf{5)} $A$ is of finite representation type and the right Auslander ring of $A$ is a generalized right QF-2 ring.

\textbf{6)} $A$ is of finite representation type and the right Auslander ring of $A$ is a generalized left co-QF-2 ring.

\textbf{7)} $A$ is of finite representation type with basic set of primitive idempotents $e_1,\ldots,e_n$ and the socle of each indecomposable right $A$-module $U$ is isomorphically embedded in $(e_1A/e_1J)\oplus\ldots\oplus (e_nA/e_nJ)$.

\textbf{Theorem 3.5 \cite{AsgBK23b}.}\\
Let all maximal right ideals of the ring $A$ are ideals. The following conditions are equivalent.

\textbf{1)} $A$ is a right co-K\"othe ring.

\textbf{2)} $A$ is a strongly right co-K\"othe ring.

\textbf{3)} $A$ is of finite representation type and every indecomposable module has a square-free socle.

\textbf{4)} $A$ is of finite representation type and every indecomposable module has a cyclic socle.

\textbf{5)} $A$ is of finite representation type and the right Auslander ring of $A$ is a generalized right QF-2 ring.

\textbf{6)} $A$ is of finite representation type and the right Auslander ring of $A$ is a generalized left co-QF-2 ring.

\textbf{Theorem 3.6 \cite{AsgBK23b}.}\\
Let $A$ is a finite dimensional algebra over a field. If $A$ is a strongly right co-K\"othe ring then $A$ is a left K\"othe ring.

A module $M$ is called an \textsf{extending} module if every submodule is essential in a direct summand of $M$.

\textbf{Theorem 3.7 \cite{AsgBK23b}.}\\
The following conditions are equivalent for a ring $A$.

\textbf{1)} $A$ is a very strongly right co-K\"othe ring.

\textbf{2)} Every right $A$-module is a direct sum of co-cyclic modules.

\textbf{3)} Every non-zero right $A$-module is a direct sum of modules with non-zero top and simple socle.

\textbf{4)} $A$ is of finite representation type and every (finitely generated) indecomposable right $A$-module has a simple socle.

\textbf{5)} Every right $A$-module is a direct sum of extending modules.

\textbf{6)} Every right $A$-module is a direct sum of uniform modules.

\textbf{7)} $A$ is of finite representation type and the right Auslander ring of $A$ is a right QF-2 ring. 

\textbf{8)} $A$ is of finite representation type and the right Auslander ring of $A$ is a left co-QF-2 ring.

If these assertions \textbf{1)}--\textbf{8)} hold, then $A$ is an Artinian, right serial ring.

A module $M$ is called \textsf{lifting} if for every submodule $N$ of $M$, there exists a direct sum decomposition $M=M_1\oplus M_2$ such that $M_1\subseteq N$ and $N\cap M_2$ is superfluous in $M2$.

\textbf{Theorem 3.8 \cite{AsgBK23b}.}\\
The following statements are equivalent for any ring $A$.

\textbf{1)} $A$ is a very strongly co-K\"othe ring.

\textbf{2)} $A$ is a very strongly K\"othe ring.

\textbf{3)} $A$ is an Artinian serial ring.

\textbf{4)} Every left and right $A$-module is a direct sum of uniform modules.

\textbf{5)} Every left and right $A$-module is a direct sum of extending modules.

\textbf{6)} Every left and right $A$-module is a direct sum of lifting modules.

\textbf{7)} $A$ is of finite representation type and the left (right) Auslander ring of $A$ is a QF-2 ring.

\textbf{8)} $A$ is of finite representation type and the left (right) Auslander ring of $A$ is a co-QF-2 ring.

\textbf{9)} Every left and right $A$-module is a direct sum of finitely generated modules with square-free top.

\textbf{Open question 3.9.} Solve the K\"othe's problem in general case.

\textbf{Open question 3.10.} Let over ring $A$ all right modules semidistributive. Is it true that all left $A$-modules are semidistributive?

\end{document}